\documentclass{edpa}
\usepackage{amsmath}
\usepackage{amssymb,amsfonts,amsthm}
\newcommand{\TryPackage}[3]{\IfFileExists{#1.sty}
{\usepackage{#1}#2}{#3}}\TryPackage{mathrsfs}
{\renewcommand{\mathcal}{\mathscr}}
{%% else try euler fonts
\TryPackage{eucal}{}{}}
\newcommand{\BIBand}[1]{and}

\theoremstyle{plain}

\newtheorem{theorem}{Theorem}[section]
\newtheorem{prop}[theorem]{Proposition}
\newtheorem{lemma}[theorem]{Lemma}
\newtheorem{cor}[theorem]{Corollary}
\newtheorem{conjecture}[theorem]{Conjecture}

\theoremstyle{definition}
\newtheorem{dfn}[theorem]{Definition}

\theoremstyle{definition} % i.e. same as plain, but theorembodyfont is rm

\newtheorem{example}[theorem]{Example}
\newtheorem{remark}[theorem]{Remark}

\newcommand{\exendproof}{\renewcommand{\qed}{\relax}\end{proof}}

\numberwithin{equation}{section}
\newcommand{\plref}[1]{{\normalfont \ref{#1}}}

\newcommand\gl{\lambda}

\newcommand\eps{\varepsilon}

\newcommand{\N}{\mathbb{N}}
\newcommand{\C}{\mathbb{C}}
\newcommand{\R}{\mathbb{R}}

\newcommand{\Z}{\mathbb{Z}}

\newcommand\cd{\mathcal{D}}
\newcommand\ce{\mathcal{E}}

\newcommand\ch{\mathcal{H}}

\newcommand\cl{\mathcal{L}}

\newcommand\cn{\mathcal{N}}

\newcommand\cs{\mathcal{S}}

\renewcommand\Im{\operatorname{Im}}

\newcommand\loc{\operatorname{loc}}

\newcommand\supp{\operatorname{supp}}

\newcommand\vol{\operatorname{vol}}

\newcommand{\cinfz}[1]{C_0^\infty(#1)}%C-unendlich Schnitte mit kompaktem
\newcommand{\cinf}[1]{C^\infty(#1)}%C-unendlich(Param)

\newcommand\DST{\displaystyle}

\newcommand{\ovl}[1]{\overline{#1}}
\newcommand{\comment}[1]{\relax}
\renewcommand{\tilde}{\widetilde}

\newcommand{\scalar}[2]{\langle #1,#2\rangle}
\newcommand{\setdef}[2]{\{ #1 \,|\, #2\}}
\newcommand{\bigsetdef}[2]{\bigl\{ #1 \,\bigm|\, #2\bigr\}}

\newcommand{\Mat}{\operatorname{M}}
\newcommand{\AC}{\operatorname{AC}}
\newcommand{\comp}{\operatorname{comp}}

\newenvironment{thmenum}{

\begin{enumerate}}{\end{enumerate}}

\newlength{\boxwidth}

\newcommand{\cLHc}{\cl^2_{\ch,\comp}(I)}
\newcommand{\LH}{L^2_\ch(I)}
\newcommand{\cLH}{\cl^2_\ch(I)}

\newcommand{\cLHcv}[1]{\cl^2_{\ch,\comp}(#1)}
\newcommand{\LHv}[1]{L^2_\ch(#1)}
\newcommand{\cLHv}[1]{\cl^2_\ch(#1)}

\newcommand{\fsys}{Y}      %notation for fundamental system
\newcommand{\Smin}{S_{\min}} % minimal extension of the slr
\newcommand{\Smax}{S_{\max}} % maximal extension of the slr
\newcommand{\csmin}{\cs_{\min}}
\newcommand{\csmax}{\cs_{\max}}
\newcommand\restr{\restriction}

\newcommand{\Dmax}{D_{\max}}
\newcommand{\Dmin}{D_{\min}}
\newcommand{\End}{\operatorname{End}}

\author{Matthias Lesch} 
\title{Essential self--adjointness of symmetric linear relations
associated to first order systems}
\address{The University of Arizona, Department of Mathematics,
617 N. Santa Rita, Tucson, AZ 85721--0089, USA}
\email{lesch@math.arizona.edu}
\abstract{The purpose of this note is to present several
criteria for essential self--adjointness. The method
is based on ideas due to Shubin.

This note is divided into two parts. The first part deals
with symmetric first order systems on the line in the
most general setting. Such a
symmetric first order system of differential equations
gives rise naturally to a symmetric linear relation in a Hilbert space.
In this case even regularity is nontrivial. We will announce
a regularity result and discuss criteria for
essential self--adjointness of such systems. 
A byproduct of the regularity result is a short proof
of a result due to Kogan and Rofe--Beketov \cite{KogRof:SIS}: the so--called formal
deficiency indices of a symmetric first order
system are locally constant on $\C\setminus\R$.
The regularity and its corollary are based on joint work
with Mark Malamud. Details will be published elsewhere.

In the second part we consider a complete Riemannian manifold,
$M$,  and a first order differential operator, $D:\cinfz{E}\to \cinfz{F}$,
acting between sections of the hermitian vector bundles $E,F$.
Moreover, let $V:\cinf{E}\to L^{\infty}_{\loc}(E)$ 
be a self--adjoint zero order differential
operator. We give a sufficient condition for the Schr\"odinger
operator $H=D^tD+V$ to be essentially self--adjoint.
This generalizes recent work of I. Oleinik \cite{Ole:ESA,Ole:CCQ,Ole:ESAG}, M. Shubin
\cite{Shu:CQC,Shu:ESA}, and M. Braverman \cite{Bra:SAS}.

We essentially use the method of Shubin. Our presentation shows
that there is a close link between Shubin's self--adjointness
condition for the Schr\"odinger operator and Chernoff's self--adjointness
condition for powers of first order operators.

We also discuss non--elliptic operators. However, in this case
we need an additional assumption. We conjecture that the additional
assumption turns out to be obsolete in general.

The criteria we are going to present in the first and 
second part of this note are very closely related. In fact, 
after we had  done the second part, we saw that the theory 
can be extended to symmetric linear relations associated 
to symmetric first order systems.
}
\msc{Primary 34L05; Secondary 35P05, 58G25}
\keyword{linear relation, self--adjoint}
\http{www.math.arizona.edu/$\sim$lesch}
\begin{document}
\maketitle
\tableofcontents

\section{First order systems on the line}

Let $I\subset\R$ be an interval. We consider a first order
system
\begin{equation}
    J(x)\frac{df}{dx}+B(x) f(x)=\ch(x) g(x),
\label{SDM-G1.1}
\end{equation}
where
\begin{align}
&J\in\AC(I,\Mat(n,\C)), &&J(x)=-J(x)^*,\quad \det J(x)\not=0, \;
\textrm{for}\; x\in I,\nonumber\\
       &B\in L^1_{\loc}(I,\Mat(n,\C)), &&
B(x)=B(x)^*-J'(x), \;\textrm{for}\; x\in
       I,\label{SDM-G1.2}\\
       &\ch\in L^1_{\loc}(I,\Mat(n,\C)),&& 
\ch(x)=\ch(x)^*, \quad \ch(x)\ge 0,\;\textrm{for}\; x\in
       I.\nonumber
\end{align}
Here, $\Mat(n,\C)$ denotes the set of complex $n\times n$ matrices and
$\AC(I,\Mat(n,\C))$ the set of absolute continuous functions with
values in $\Mat(n,\C)$.

We need some more notation: we equip $C_0(I,\C^n)$, the space
of continuous $\C^n$--valued functions with compact support, with
the (semidefinite) scalar product
   \begin{equation}
\scalar{f}{g}_\ch:=\int_I f(x)^*\ch(x)g(x)dx,
\label{SDM-G1.3}
    \end{equation}
and denote by $\cLH$ the completion of $C_0(I,\C^n)$
with respect to the semi-norm induced by \eqref{SDM-G1.3}.
Alternatively, $\cLH$ can be described as the set of Borel--measurable
$\C^n$--valued functions satisfying
$\scalar{f}{f}_\ch:=\int_I f(x)^*\ch(x)f(x)dx <\infty.$
As usual, one puts 
%\begin{equation}
$\LH:=\cLH/\setdef{f\in\cLH}{\|f
\|_\ch=0}$.
%\end{equation}
$\LH$ is a Hilbert space. For a function $f\in \cLH$ 
we will
denote by $\tilde f$ the corresponding class in $\LH$. 
If $\ch(x)$ is invertible a.e. then a class $\tilde f$ 
contains
at most one continuous representative, hence if 
$\ch(x)$ is invertible a.e.
and $f$ is continuous we 
will not distinguish between $f$ and $\tilde f$.

Assume for the moment that $\ch(x)$ is invertible 
for almost all $x\in I$ and $\ch(x)^{-1}\in 
L^1_{\loc}(I,\Mat(n,\C))$. Then \eqref{SDM-G1.1} induces a 
symmetric operator
   \begin{equation}
    L:=\ch^{-1}\Big(J \frac{d}{dx}+B\Big)
\label{SDM-G1.4}
  \end{equation}
in the Hilbert space $\LH$ with domain 
$\cd(L)=\AC_{\comp}(I,\C^n)$. The symmetry is 
implied by $B=B^*-J'$ and $\ch^*=\ch$.
However, the interesting case is the one where $\ch(x)$ is singular.
If $\ch(x)$ is singular then \eqref{SDM-G1.1} will in general neither
define an operator nor will it be densely defined. Rather it 
will  give rise to a
\emph{symmetric linear relation}, whose definition we
recall for the reader's convenience:
\begin{dfn}\label{SDM-S1.1} Let $\mathfrak H$ be a linear space equipped with a positive
semidefinite hermitian sesqui--linear form $\scalar{\cdot}{\cdot}$.
A linear subspace $\cs\subset \mathfrak{H}\times\mathfrak{H}$ is called
a symmetric linear relation (s.l.r.) if for $\{f_j,g_j\}\in \cs, j=1,2$, one has
$\scalar{f_1}{g_2}=\scalar{f_2}{g_1}$.
\end{dfn}

For example, the graph of an (unbounded) symmetric operator in $\mathfrak H$ is
a s.l.r. 
The system \eqref{SDM-G1.1} defines a symmetric linear relation, 
$\csmin$, in  
$\cLH$ as
follows: $\{f, g\}\in \csmin$ if and only if $f\in
\AC_{\comp}(I,\C^n), g\in \cLHc$ and
$Jf'+Bf=\ch g$.

$\csmin$ induces a symmetric linear relation, $\Smin$, in 
$\LH$ in a fairly
straightforward way: $\{\tilde f,\tilde g\}\in \Smin$ if 
and only if
there exist representatives $f\in \tilde f, g\in 
\tilde g$ such that
$\{f,g\}\in \csmin$. 

Looking at first order systems seems to be rather special. Therefore, it
is important to note that 
an arbitrary symmetric $n^{th}$--order system is unitarily
equivalent to a symmetric first order system (\cite{KogRof:SIS},
\cite{Orc:CDE}). 
In most cases, however, the Hamiltonian $\ch$ of this first order system will be singular.
As an example we show how a second order Sturm--Liouville equation
can be transformed into a system of the form \eqref{SDM-G1.1}:

\begin{example}\label{SDM-S1.2}
We consider a Sturm--Liouville type equation
   \begin{equation}
   -\frac{d}{dx}\bigl(A(x)^{-1}\frac{d}{dx} u(x)\bigr)+V(x)u(x)=\ch(x) v(x),
\label{SDM-G1.5}
    \end{equation}
where $A,V,\ch\in L^1_{\loc}(I,\Mat(n,\C))$ and $A(x)$ is positive definite
for all $x\in I$. The system \eqref{SDM-G1.5} defines a symmetric linear
relation as follows: $\{u,v\}\in\csmin$ if and only if $u\in\AC_{\comp}(I,\C^n),
A^{-1}\frac{d}{dx}u\in \AC_{\comp}(I,\C^n),$
$v\in \cLHc$ and \eqref{SDM-G1.5} holds. As before, let
$\Smin:=\setdef{\{\tilde u,\tilde v\}}{\{u,v\}\in\csmin}$.
%in the distributional sense.
%Note that this implies
%that $A^{-1}\frac{d}{dx}u\in \AC_{\comp}(I,\C^n)$
%(if $v\in\cLHc$ then, since $\ch\in L^1_{\loc}(I)$,

% $\ch v\in L^1_{\comp}(I)$).
Note that if $v\in\cLHc$ then, since $\ch\in L^1_{\loc}(I,\Mat(n,\C))$,
$\ch v\in L^1_{\comp}(I,\C^n)$.
Consequently,
$\{(u,iA^{-1}\frac{d}{dx}u),(v,0)\}$ is in the symmetric linear
relation, $\tilde\csmin$, induced by the system
  \begin{equation}
\begin{pmatrix} 0 & i\\i&0\end{pmatrix}
\begin{pmatrix} f_1\\f_2\end{pmatrix}'
+\begin{pmatrix} V&0\\0&-A\end{pmatrix}
\begin{pmatrix} f_1\\f_2\end{pmatrix}
=
\begin{pmatrix} \ch&0\\0&0\end{pmatrix}
\begin{pmatrix} g_1\\g_2\end{pmatrix}.
\label{SDM-G1.6}
   \end{equation}
Conversely, if $\{(f_1,f_2),(g_1,g_2)\}\in\tilde\csmin$ then
$\{f_1,g_1\}\in\csmin$.
It is also clear that the Hilbert spaces $\LH$ and $L^2_{\tilde\ch}(I)$,
$\tilde\ch=\begin{pmatrix} \ch&0\\0&0\end{pmatrix}$, are canonically
isomorphic. Hence the s.l.r. $\Smin$ and $\tilde S_{\min}$
in $\LH$ resp. $L^2_{\tilde\ch}(I)$ are unitarily equivalent.

If $\ch(x)$ is invertible and $\ch(x)^{-1}\in L^1_{\loc}(I,\Mat(n,\C))$ then
$\Smin$ is (the graph of) a densely defined symmetric operator in the 
Hilbert space $\LH$. However, $\tilde\ch(x)$ is singular everywhere.
\end{example}

The following example shows that the domain of the s.l.r. $\Smin$ can be rather
small:

\begin{example}\label{SDM-S1.3} Let $B=0, J=\begin{pmatrix} 0&1\\ 
-1& 0
\end{pmatrix},$ and $ \ch(x)=\begin{pmatrix} 1 & 0\\ 
0&0\end{pmatrix}$.
If $\{f,g\}\in \csmin$ then $f_2'=g_1, f_1'=0$, and since
$f$ is absolute continuous with compact support we 
infer $f_1=0$. Hence $\ch f=0$ and thus $\tilde 
f=0$. Thus, the domain of $\Smin$ is $\{0\}$.
\end{example}

%\begin{dfn} We denote by $\scmin(J,B,\ch)$ the s.l.r. associated to the system
%\eqref{SDM-G1.1} in case the dependence on $J,B,\ch$ has to be emphasized.
%Furthermore, let $

The system \eqref{SDM-G1.1} can be simplified further and put into canonical form.
Details of the construction can be found in \cite[Sec. 1.3]{KogRof:SIS}
or \cite{LesMal:NSS}. For the moment denote by 
$\cs(J,B,\ch)$ the s.l.r. induced by the system \eqref{SDM-G1.1}. 
A "gauge  transformation" $U\in \AC(I,\operatorname{GL}(n,\C))$ 
induces a unitary map 
\begin{equation}%\begin{split}
\Psi_U:\cLH\to \cl_{\tilde \ch}(I), \quad f\mapsto 
U^{-1} f,\quad \tilde \ch:= U^*\ch U,
\label{SDM-G1.7}
\end{equation}
and a simple computation shows that 
\begin{equation}
    \Psi_U \cs(J,B,\ch)\Psi_U^*= \cs(\tilde J,\tilde 
B,\tilde \ch),
\label{SDM-G2.8}
\end{equation}
where
\begin{equation}\label{SDM-G2.9}%\begin{split}
    \tilde J= U^*JU,\quad \tilde B= U^*JU'+U^*BU,\quad \tilde \ch= U^*\ch U.
%			       \end{split}
\end{equation}
It can be shown that the gauge transformation $U$ can be chosen
in such a way that $J$ is constant and $B=0$. Such a system is called
"canonical". 

Pick $x_0\in I$ and let 
$\fsys(.,\gl):I\rightarrow \Mat(n,\C)$ be the solution 
of the initial value problem
\begin{equation}
     J(x)\fsys'(x,\gl)+B(x)\fsys(x,\gl)=\gl \ch(x) 
\fsys(x,\gl),\quad \fsys(x_0,\gl)=I_n.
\end{equation}
Here, $I_n$ denotes the $n\times n$ unit matrix. The 
existence of $\fsys$ follows from the integrability assumptions
in \eqref{SDM-G1.2}.

\begin{dfn}\label{SDM-S1.4} The system \eqref{SDM-G1.1} is 
said to be
\emph{definite} on $I$ if there exists a compact subinterval
$I_0\subset I$ such that the matrix
\begin{equation}
\int_{I_0} \fsys(x,\gl)^*\ch(x) \fsys(x,\gl)dx
\label{SDM-G1.11}
\end{equation}
is invertible for a $\gl\in\C$.
\end{dfn}
If the system is definite then \eqref{SDM-G1.11} is invertible
for all $\gl\in\C$ \cite[Theorem 1.1]{KogRof:SIS}.
The property of a system \eqref{SDM-G1.1} to be definite 
is gauge invariant. There is a simple criterion for definiteness:
namely, if there exists a compact subinterval $I_0\subset I$
such that $\int_{I_0}\ch$ is invertible, then the system is
definite. For a canonical system ($B=0$) this criterion is
also necessary. In general, the definiteness
will also depend on $J$ and $B$.

%%%%%%%%%%%%%%%%%%%%%%%%%%%%%%%%%%%%%%
Some bibliographic comments are in order, however we do not claim
to give a complete historical account: 
A standard reference for symmetric linear relations arising from
symmetric first order systems is the thesis of Orcutt \cite{Orc:CDE}, 
which unfortunately has not been published. 
Other references are \cite{Ben:SRH}, \cite{LanTex:GKM},
\cite{DijSno:SER}. First order systems have been studied
extensively in \cite{KogRof:SIS}. Canonical systems are discussed 
in great detail in \cite{GohKre:TAV}.
%%%%%%%%%%%%%%%%%%%%%%%%%%%%%%%%%%%%%%%%%%%%%%%%%%%%%%%%%%%%%%%%%%%%%%%%%%%%

\subsection{Regularity of the maximal relation}

We consider again the system \eqref{SDM-G1.1}, \eqref{SDM-G1.2}.

\begin{dfn}\label{SDM-S1.5} We denote by $S$ the closure in $\LH\times \LH$
of $\Smin$ and by $\Smax:=S^*=\setdef{\{f,g\}\in \LH\times\LH}{
\scalar{f}{v}=\scalar{g}{u}\;\text{for all}\; \{u,v\}\in S}$ 
the adjoint of $S$. Moreover, let 
\[\csmax:=\bigsetdef{\{f,g\}}{f,g\in\cLH, f\in\AC(I,\C^n), Jf'+Bf=\ch g}.\]
\end{dfn}

The notation $\csmax$ is deliberately chosen:
if $S$ is the graph of a symmetric first order operator
as in \eqref{SDM-G1.4} then it is well--known that 
each pair $\{\tilde f,\tilde g\}$ has representatives
$\{f,g\}\in\csmax$. It is exaggerating but true that
this follows from elliptic regularity. 
For the system \eqref{SDM-G1.1} the same statement holds true,
although it is less obvious:

\begin{theorem}[Regularity Theorem] \label{SDM-S1.6} Let $\{\tilde f,\tilde g\}\in\Smax$.
Then for each representative $g\in\tilde g$ there exists $f\in\tilde f$
such that $\{f,g\}\in\csmax$.
\end{theorem}

For definite systems 
this has been proved by Orcutt \cite[Thm. II.2.6 and 
Thm. IV.2.5]{Orc:CDE}.
Another proof for (not necessarily definite) $2\times 
2$ canonical systems was given by
I.S. Kac \cite{Kac:LRGC} in the deposited but unpublished
elaboration of \cite{Kac:LRG}. 
The proof of a more detailed version of Theorem \plref{SDM-S1.6}
will be published in \cite[Sec. 2]{LesMal:NSS}.

We present an application of the regularity theorem: Let
\begin{equation}
      \ce_{\gl}(S):=\bigsetdef{f\in\cLHv{I}\cap 
AC(I,\C^n)}{ Jf'+ 
Bf=\gl \ch f},
       \label{SDM-G1.12}
\end{equation}
and denote by
%\begin{equation}
 $    \cn_{\pm}(S):=\dim \ce_{\pm i}(S)$
%\end{equation}
the \emph{formal} deficiency indices of the system 
\eqref{SDM-G1.1}. 
Furthermore, for a symmetric linear relation $A$ in 
the Hilbert space $\mathfrak H$ we denote by
   \begin{equation}
        E_\gl(A):=\bigsetdef{f\in \mathfrak 
H}{\{f,\gl f\}\in A^*},\qquad \gl\in\C,
\label{SDM-G1.13}
   \end{equation}
the defect subspace and by
%   \begin{equation}
 $       N_\pm(A):=\dim E_{\pm i}(A)$
%    \end{equation}
the deficiency indices of $A$. 
It is well--known  
that  
\begin{equation}
\dim E_{\pm \gl}(A)=N_{\pm}(A),\qquad \gl \in \C_+:=\bigsetdef{z\in\Z}{\Im
z>0}.
\label{SDM-G1.14}
\end{equation}
Namely, the relation $A^*-\gl$ is semi--Fredholm
for $\gl\in\C\setminus\R$. Thus $\dim E_\gl(A)$ is locally 
constant on $\C\setminus \R$ and therefore 
$\dim E_{\pm\gl}(A)=\dim E_{\pm i}(A)$ 
for $\gl\in \C_+$. 

The same statement for the dimensions of the formal 
defect subspaces $\ce_{\gl}(S)$ is true but less trivial.
The only proof we know of is due to Kogan and
Rofe--Beketov \cite[Sec. 2]{KogRof:SIS}. It
uses methods from complex analysis and is rather technical.
%First of them follows from the observation that the 
%relation $A^*-\lambda$
%is semi--Fredholm for $\lambda \in {\Bbb C}\setminus 
%{\Bbb R}.$ This fact 
%yields (see [Kato]) locally constancy of $dim{\mathfrak 
%N}_{\lambda }(A)$ on
%${\Bbb C}\setminus {\Bbb R}$ which in turn yields the 
%required equalities.
%Since for $\gl\in \C\setminus\R$ the relation 
%$A^*-\gl$ is semi--Fredholm, 
%it is clear that $\dim {\mathfrak N}_\gl(A)$ is locally 
%constant on 
%$\C\setminus \R$ and therefore 
%dim${\mathfrak N}_{\pm\gl}(A)=$dim${\mathfrak N}_{\pm i}(A)$ 
%for $\gl\in {\Bbb C}_+.$ 
Using Theorem \plref{SDM-S1.6} we can give a painless proof of this
fact:
   
\begin{theorem}[{\cite[Theorem 
2.1]{KogRof:SIS}, \cite[Sec. 2]{LesMal:NSS}}]\label{S1.11}
%\marginpar{exact reference}
Let $S$ be a general symmetric system \eqref{SDM-G1.1}, \eqref{SDM-G1.2} on
an interval $I\subset \R$. 
If the system is definite or if the interval is half--closed, i.e. $I=[0,a)$,
then 
\[\dim \ce_{\pm \lambda }(S)=\dim \ce_{\pm i}(S)=:\cn_{\pm 
}(S),\quad\text{for}\quad \gl\in\C_+.
\]
\end{theorem}
%%%%%%%%%%%%%%%%%%%%%
\begin{proof}
1. We assume first that the system $S$ is definite. Then 
the quotient map $\ce_\gl(S)\to E_\gl(S), f\mapsto \tilde f$
is bijective. 

Indeed, the injectivity follows immediately from the
definition of definiteness. To prove surjectivity, 
consider $\tilde f\in E_\gl(S)$. This means $\{\tilde f,\gl\tilde f\}\in \Smax$
and in view of Theorem \plref{SDM-S1.6} there exists
$f\in\tilde f,f\in \AC(I,\C^n)\cap \cLH$ 
such that $Jf'+Bf=\gl\ch f$.
Thus $f\in E_{\lambda }(S)$. 
This proves surjectivity.

Now we have $\dim \ce_{\gl}(S)=\dim E_{\gl}(S)$ 
and in view of \eqref{SDM-G1.14} we reach the conclusion. 

2. If $S$ is not definite but $I=[0,a)$ we replace 
${\ch}$ by $\tilde \ch={\ch}+\chi I_n$, where $\chi $ is the 
characteristic function of an interval $[0,\eps)\subset I$. 
The system $\tilde S=S(J,B,\tilde \ch)$ is definite on 
$I$ and 1. applies. To complete the proof it remains
to note that we obtain a linear
isomorphism, $\Phi$, from $\ce_\gl(S)$ onto $\ce_\gl(\tilde S)$ as follows:
for $f\in \ce_\gl(S)$ let $\Phi f$ be the solution of the differential
equation $Jy'+B y=\gl \tilde \ch y$ with $\Phi f\restr [\eps,a)=f\restr [\eps,a)$.
\end{proof}

\subsection{Essential self--adjointness}

In this section we study the system \eqref{SDM-G1.1} on 
the real line and discuss criteria for essential self--adjointness. 
%\subsection{Essential self--adjointness in the case of
%a singular Hamiltonian}
%Next we want to present a criterion for essential 
%self--adjointness
%in a case where the Hamiltonian is singular 
%everywhere. 
%In particular we will show that in the situation of 
%Theorem \plref{S2.2} the square of $\cs$ is 
%essentially self--adjoint,
%too. Actually, the method presented here is very close
%to recent work of Shubin \cite{Shu:CQC}, see also
%Lesch \cite{Les:RPS} for a generalization.
As a motivation, let $\{f,h\}$ be in the "square"
of $\csmin$, that is there is a $g\in\cLH$ such that 
$\{f,g\}\in\csmin$ and $\{g,h\}\in\csmin$. This is 
equivalent to the equation
\begin{equation}
\begin{pmatrix} 0 & J\\ J & 0\end{pmatrix}
\begin{pmatrix} f \\ g \end{pmatrix}'
+\begin{pmatrix} 0 & B \\ B & -\ch \end{pmatrix}
\begin{pmatrix} f\\ g\end{pmatrix}
= 
\begin{pmatrix} \ch & 0\\0 & 0\end{pmatrix}
\begin{pmatrix} h \\ 0\end{pmatrix},
\label{SDM-G1.15}
   \end{equation}
with $f,g\in\AC_{\comp}(I,\C^n), h\in\cLHc$.
A second example is the system discussed in Example
\plref{SDM-S1.2}. These examples lead
us to consider a first order system
   \begin{equation}
\tilde J f'+\tilde B f=\tilde \ch g,
\label{SDM-G1.16}
  \end{equation}
where
  \begin{equation}
\tilde J=\begin{pmatrix} 0 & J\\ J & 0\end{pmatrix}
,\quad
\tilde B=\begin{pmatrix} V & B \\ B & -A \end{pmatrix},
\quad
\tilde\ch=\begin{pmatrix} \ch & 0\\0 & 0\end{pmatrix}.
\label{SDM-G1.17}
  \end{equation}
$A$ is assumed to be nonnegative. $V$ may be viewed as a "potential"
added to $\csmin^2$. It is clear that $L^2_{\tilde \ch}(I)$ is canonically
isomorphic to $L^2_{\ch}(I)$. We put 
$\tilde\csmin=\cs(\tilde J,\tilde B,\tilde \ch)$.
For simplicity we will consider the interval $\R$ only.
For a function $f\in \cl_{\tilde \ch}^2(\R)$ we denote
by $f_1,f_2$ the first resp. last $n$ components.
%We denote by $r(x)$ the function defined in \eqref{ML-G3.2} with respect to

%$J$ and $\ch$.
We will use several times that if $\ch(x)$ and $A(x)$ are invertible
then we can estimate, for $\xi,\eta\in\C^n$,
   \begin{equation}
\begin{split}
   \big| \xi^* J\eta|&=\|A(x)^{1/2}\xi\|
   \|A(x)^{-1/2}J(x)\ch(x)^{-1/2}\ch(x)^{1/2}\eta\|\\
    &\le  \|A(x)^{-1/2}J(x)\ch(x)^{-1/2}\| \|A(x)^{1/2}\xi\| \|\ch(x)^{1/2}\eta\|.
\end{split}
\label{SDM-G1.18}
   \end{equation}
Thus we put
   \begin{equation}
    c(x):=\begin{cases}  \|A(x)^{-1/2}J(x)\ch(x)^{-1/2}\|,&
    \det(A(x)\ch(x))\not=0,\\
          \infty,&\textup{otherwise}.
          \end{cases}
\label{SDM-G1.19}
     \end{equation}
The self--adjointness criterion we are going to
present will depend also on $V$. We assume
that there exists an absolute continuous function $q\ge 1$
on $\R$ such that
    \begin{equation}
    V\ge -q \ch.
\label{SDM-G1.20}
   \end{equation}
%%%%%%%%%%%%%%%%%%%%%%%%%%%%%%%%%%

\begin{lemma}\label{SDM-S1.8}
Let $f\in L^1_{\loc}(\R), f(x)\ge 0$, be a non--negative
locally integrable function.
Assume in addition that
\begin{equation}
\pm\int_0^{\pm\infty} f(x) dx=+\infty.
\label{SDM-G1.21}
\end{equation}
Then there is a sequence of functions $\chi_n\in\AC_{\comp}(\R)$
satisfying
\begin{equation}%\begin{split}
         0\le \chi_n\le 1,\quad |\chi_n'|\le \frac 1n f(x),\quad
         \lim_{n\to\infty} \chi_n(x)=1,\quad x\in \R.
		%\end{split}
\label{SDM-G1.22}
\end{equation}
\end{lemma}
\begin{proof}
Let $\chi\in\cinfz{\R}$ with $0\le\chi\le 1, $
$\chi(x)=1$ in a neighborhood of $0$ and $|\chi'|\le 1$.
Then
\begin{equation}
   \chi_n(x):=\chi\Bigl( \frac 1n \int_0^x f(s)ds\Bigr)
\end{equation}
does the job.
\end{proof}
\comment{Fix $n\in\N$. By B. Levy's theorem on monotone 
convergence we have
\[ \lim_{C\to+\infty}\int_0^\infty \frac 1n 
\min(C,f(x))dx=+\infty,\]
and thus we may choose $C>0$ such that 
\[   \int_n^\infty  \min(C,\frac 1n f(x))dx \ge 2.\]
Now choose $N$ large enough such that
\[ K_n:=\int_n^N  \min(C,\frac 1n f(x)) dx \ge 1\]
and put
\[\chi_n(x):=1- \frac {1}{K_n} 
\int_{\min(n,x)}^{\min(N,x)}\min(C,\frac 1n f(s))ds.\]
$\chi_n$ has the desired properties with $x_n=N$.
   \end{proof}}
%%%%%%%%%%%%%%%%%%%%%%%%%%%%%%%%

\begin{lemma}[{\cite[Lemma 3.1]{Shu:CQC}, cf. Proposition \ref{S5} below}]
\label{SDM-S1.9} Assume that
\[ \pm\int_0^{\pm\infty} \frac{1}{c(x)}dx=\infty,\]
and that $|\frac{d}{dx}q^{-1/2}(x)|\le C/c(x).$
Let $\{f,g\}\in\tilde\cs_{\max}$. 
Then
$q^{-1/2}f_2\in\cl_{A}^2(\R)$ and
\[  \|q^{-1/2}f_2\|_{A}\le 2\Bigl((1+2C^2)\|f\|_{\tilde \ch}^2+
\|f\|_{\tilde \ch}\|g\|_{\tilde \ch}\Bigr).\]
\end{lemma}
\begin{proof} By Lemma \plref{SDM-S1.8} there are absolute continuous functions
$\chi_n$  with $0\le \chi_n\le 1$, $\lim\limits_{n\to\infty}\chi_n(x)=1$, and
\begin{equation}   |\chi_n'(x)| \le \frac {1}{nc(x)}.\end{equation}
Put $\psi_n:=\chi_n q^{-1/2}$. We have
\begin{equation}
  |\psi_n'(x)|\le \bigl(\frac 1n +C\bigr)\frac{1}{c(x)}=: C_n\frac{1}{c(x)}.
\label{SDM-G1.24}
\end{equation}
Then
% \begin{equation}
% \begin{split}
\begin{align}
    \|&\psi_n f_2\|_{A}^2= \int_{\R} \psi_n^2(x)
      f_2^*(x)(Jf_1'+Bf_1)(x)dx\nonumber \\
      &=\int_{\R}\psi_n^2 (J(x) f_2'(x)+B(x)f_2(x))^*f_1(x)-2
         \int_{\R} \psi_n(x)\psi_n'(x) f_2(x)^* J(x) f_1(x)dx \nonumber \\
      &=\int_{\R} \psi_n(x)^2 g_1(x)^*\ch(x) f_1(x) dx
          -\int_{\R} \psi_n(x)^2 f_1(x)^*V(x)f_1(x) dx \nonumber \\
   &\quad - 2 \int_{\R}
         \psi_n(x) \psi_n'(x)f_2(x)^* J(x) f_1(x)dx.
% \end{split}
\label{SDM-G1.25}
\end{align}
%     \end{equation}
Note that in view of
\eqref{SDM-G1.24}  the matrices $A(x)$ and $\ch(x)$ are invertible if
$\psi_n'(x)\not=0$. Combining \eqref{SDM-G1.18}, \eqref{SDM-G1.20},
\eqref{SDM-G1.24}, \eqref{SDM-G1.25} and the well--known
estimate $2|ab|\le a^2+b^2$ we obtain
\begin{equation}\begin{split}
 \|\psi_n f_2\|_{A}^2&\le \big|\scalar{\psi_n^2 f}{g}_{\tilde\ch}\big|
                +\|\psi_n q^{1/2} f_1\|_{\ch}^2
                +2C_n\|\psi_n f_2\|_{A}\|f\|_{\ch}\\
        &\le \|f\|_{\tilde\ch}\|g\|_{\tilde\ch}+
           (1+2C_n^2)\|f\|_{\tilde\ch}^2
                   +\frac 12 \|\psi_n f_2\|_{A}^2,
                \end{split}
\label{SDM-G1.26}
\end{equation}
or
\begin{equation}
  \|\psi_n f_2\|_{A}^2
              \le  2\bigl((1+2C_n^2)\|f\|_{\tilde\ch}^2
                     +\|f\|_{\tilde\ch}\|g\|_{\tilde\ch}\bigr).
\label{SDM-G1.27}
   \end{equation}
Letting $n\to\infty$ we reach the conclusion.
  \end{proof}
%%%%%%%%%%%%%%%%%%%%%%%%%%%%%%%%%%%%
%%%%%%%%%%%%%%%%%%%%%%%%%%%%%%%%%%%%%
\begin{theorem}[{\cite[Theorem 1.1]{Shu:CQC}, cf. Theorem \ref{main} below}]\label{SDM-S1.10}
On the interval $\R$ let $\tilde J,\tilde B,\tilde \ch$ be as in \eqref{SDM-G1.17} with
$A\ge 0$.
Let $q\ge 1$ be absolute continuous and $V\ge -q\ch$.
Moreover, assume that
\begin{thmenum}
\item $\big|\frac{d}{dx}q^{-1/2}(x)\big|\le \frac{C}{c(x)}$.
\item $\displaystyle \pm \int_0^{\pm\infty} \frac{1}{c(x) q^{1/2}(x)} dx=\infty$.
\end{thmenum}
Then $\tilde S=S(\tilde J,\tilde B,\tilde \ch)$ is essentially self--adjoint.
\end{theorem}
\begin{proof}
By Lemma \plref{SDM-S1.8} there are absolute continuous functions
$\chi_n\in\AC_{\comp}(\R)$, $0\le\chi_n\le 1$,
$\lim\limits_{n\to\infty}\chi_n(x)=1$,
and
\begin{equation}
|\chi_n'(x)| \le \frac {1}{n c(x)q^{1/2}(x)}.
\end{equation}
Note that, again, $\chi_n'(x)\not=0$ implies that $A(x)$ and $\ch(x)$
are invertible. 
%Consider $\{f,g\}, \{u,v\}\in\tilde\cs_{\max}$. 
In view of the regularity
Theorem \plref{SDM-S1.6} it suffices to show for 
$\{f,g\}, \{u,v\}\in\tilde\cs_{\max}$ that
   \begin{equation}
      \scalar{f}{v}=\scalar{g}{u}.
  \end{equation}
By dominated convergence we have
\begin{equation}
\lim_{n\to\infty} \big(\scalar{\chi_n f}{v}-\scalar{\chi_n g}{u}
\big)= \scalar{f}{v}-\scalar{g}{u}.
\end{equation}
Integration by parts shows that
  \begin{equation}
\begin{split}
     \bigl(\scalar{\chi_n f}{v}-\scalar{\chi_n g}{u}\bigr)&=
       -\int_{\R}\chi_n'(x) f(x)^*\tilde J(x) u(x)dx\\
      &= -\int_{\R} \chi_n'(x)\big(f_1(x)^* J(x) u_2(x)+f_2(x)^*J(x)u_1(x)\big)dx.
\end{split}
\end{equation}
Using \eqref{SDM-G1.18} and Lemma \plref{SDM-S1.9} this
can be estimated by
  \begin{equation}
\big|\scalar{\chi_n f}{v}-\scalar{\chi_n g}{u}\big|
\le \frac 1n \big(\|f_1\|_{\ch}\|q^{-1/2}u_2\|_{A}+
                   \|q^{-1/2}f_2\|_{A}\|u_1\|_{\ch}\big),
\end{equation}
and we reach the conclusion.
 \end{proof}
%%%%%%%%%%%%%%%%%%%%%%%%%%%%%%%%%%%%%%%%%%%%%
\begin{remark}\label{SDM-S1.11}%Umschreiben in Corollar!!!!!
%1. Note that if $V\ge 0$ we may choose $q=1$. Then
%Lemma \plref{ML-S3.2.1} and Theorem \plref{ML-S3.2.2}
%hold under the only condition $\pm\int_0^{\pm\infty}\frac{1}{c(x)}dx=+\infty$.
%This shows in particular that the square of $\cs$ is essentially self--adjoint
%under the same condition as in Theorem \plref{S2.2}.
We emphasize that Lemma \plref{SDM-S1.9},
Theorem \plref{SDM-S1.10} and their proofs are adapted from
a method due to M. Shubin \cite{Shu:CQC}
who proved essential self--adjointness
for certain Schr\"odinger type operators on complete manifolds.
A generalization of Shubin's method is presented below in the
second part of this paper.
\end{remark}
We single out some special cases of the previous theorem.

\begin{cor}\label{SDM-S1.12}
Consider the system $\Smin=S(J,B,\ch)$ as in \eqref{SDM-G1.1}
on $I=\R$. Put
 \begin{equation}
    c(x):=\begin{cases}  \|\ch(x)^{-1/2}J(x)\ch(x)^{-1/2}\|,&
    \det(\ch(x))\not=0,\\
          \infty,&\textup{otherwise}.
          \end{cases}
\label{SDM-G1.34}
     \end{equation}
Assume
\begin{equation}
 \pm\int_0^{\pm\infty} \frac{1}{c(x)} dx=+\infty.
\end{equation}
Then $\Smin$ and $\Smin^2$ are essentially self--adjoint, i.e. 
$\overline{\Smin}=\Smax$ and $\overline{\Smin^2}=(S^2)_{\max}$.
\end{cor}
This corollary generalizes a result of Sakhnovich \cite{Sak:DIS}.
\begin{proof} The essential self--adjointness of $\Smin^2$
follows, in view of \eqref{SDM-G1.15}, from Theorem
\plref{SDM-S1.10} with $V=0, q=1$ and $A=\ch$.

It is easy to see that, as in the case of a symmetric
operator, the essential self--adjointness of the square of a
s.l.r. in a Hilbert space implies the essential self--adjointness of the s.l.r. itself.
However, the essential self--adjointness of $\Smin$ can easily
be seen directly:

According to Lemma \plref{SDM-S1.8} let 
$\chi_n\in\AC_{\comp}(\R)$ with $0\le\chi_n\le 1, \lim\limits_{n\to\infty}\chi_n(x)=1$,
and
\begin{equation}  |\chi_n'(x)| \le \frac {1}{n c(x)}.
\end{equation}
For $\{\tilde f,\tilde g\}\in \Smax$ we choose, 
according to Theorem \plref{SDM-S1.6}, representatives $\{f,g\}\in\csmax$ and 
put $f_n:= \chi_n f$.
Since $\chi_n'$ vanishes if $\ch(x)$ is not invertible 
the function
$\chi_n' \ch(x)^{-1} Jf$ is well--defined. Moreover
\begin{align*}
    \|\chi_n' \ch^{-1} J f\|_{\LHv{\R}}^2&\le \int_\R 
|\chi_n'(x)|^2
                f(x)^*J(x)^* \ch(x)^{-1}J(x) f(x) dx\\
     &\le \sup_{x\in \R}( \chi_n'(x)c(x))^2 
\|f\|_{\LHv{\R}}^2\\
     &\le \frac{1}{n^2} \|f\|_{\LHv{\R}}^2,
\end{align*}
hence $\chi_n' \ch(x)^{-1} Jf$ lies in $\cLHv{\R}$ and 
it converges to $0$ in $\cLHv{\R}$.
Finally, we calculate
\begin{align*}
     Jf_n' +Bf_n &= \chi_n (J f' + Bf)+ \chi_n' Jf\\
                &= \ch(\chi_n g+\chi_n' \ch^{-1} Jf)\\
                &=:\ch g_n.
\end{align*}
Thus $\{f_n,g_n\}\in \csmin$ and 
$\lim\limits_{n\to\infty}\{\tilde f_n,\tilde
g_n\}=\{\tilde f,\tilde g\}$ and the claim is proved.
\end{proof}

\begin{cor}\label{SDM-S1.13}
Let $\Smin$ be the symmetric linear relation in $L^2_{\ch}(\R)$ induced
by the Sturm--Liouville type equation
  \begin{equation}
    -\frac{d}{dx}\bigl(A(x)^{-1}\frac{d}{dx} u(x)\bigr)+V(x)u(x)=\ch(x) v(x).
\label{SDM-G1.36}
  \end{equation}
That is, $\{\tilde u,\tilde v\}\in \Smin$ if and only if there
exist $u\in\tilde u, v\in\tilde v$ such that
$u,A^{-1}\frac{d}{dx}u\in\AC_{\comp}(\R,\C^n), v\in \cLHcv{\R}$
and \eqref{SDM-G1.36} holds.
Here, we assume that $A,V,\ch\in L^1_{\loc}(\R,\Mat(n,\C))$, $A(x)$ is positive definite
for all $x\in \R$, and that there exists an absolute continuous
function $q\ge 1$ such that $V\ge -q\ch$. Let $c(x)$ be defined
by  \eqref{SDM-G1.19}.
Moreover, assume that
\begin{thmenum}
\item $\big|\frac{d}{dx}q^{-1/2}(x)\big|\le \frac{C}{c(x)}$.
\item $\displaystyle \pm \int_0^{\pm\infty} \frac{1}{c(x) q^{1/2}(x)} dx=\infty$.
\end{thmenum}
Then $\Smin$ is essentially self--adjoint.
        \end{cor}
%%%%%%%%%%%%%%%%%%%%%%%%%%%%%%%%%%%%
\begin{proof} This follows immediately from Theorem \plref{SDM-S1.10},
\eqref{SDM-G1.5}, and \eqref{SDM-G1.6}.
\end{proof}
\begin{prop}\label{SDM-S1.14} Under the assumptions of Theorem
\plref{SDM-S1.10} the
system $\tilde S=S(\tilde J,\tilde B,\tilde\ch)$ is definite.
\end{prop}
\begin{proof} Consider $f\in\cl_{\tilde\ch}^2(\R)\cap \AC(\R,\C^{2n})$
satisfying
   \begin{equation}
\tilde J f'+\tilde B f=0,\qquad \int_{\R} f^*\tilde \ch f=0.
 \label{SDM-G1.37}
  \end{equation}
We have to show that $f=0$.
\eqref{SDM-G1.37} translates into
\begin{align}
      J f_1'+ Bf_1-A f_2&=0,\label{SDM-G1.38a}\\
      J f_2'+Bf_2+V f_1  &=0,\label{SDM-G1.38b}\\
          \int_{\R} f_1^*\ch f_1&=0.\label{SDM-G1.38c}
\end{align}
Note that condition (2) in Theorem \plref{SDM-S1.10} implies that
$A(x)$ and $\ch(x)$ are invertible on a set of positive Lebesgue
measure. Consequently, the systems $\cs(J,B,A), \cs(J,B,\ch)$ are
definite.
%Consequently \eqref{ML-G3.2.15c} implies
%\begin{equation}

%          \bigl(\frac{d}{dx}q^{-1/2}\bigr) f_1=0 \;\textup{(a.e.)}.

%\label{ML-G3.2.16}

%\end{equation}

%Taking this and $V\ge -q\ch$ into account, we obtain

%\begin{equation}

%   \begin{split}

%        0&= \int_{\R} f_1^*\ch f_1\ge -\int_{\R} q^{-1} f_1^*Vf_1

%=\int_{\R} q^{-1} f_1^*(Jf_2'+Bf_2)\\

%         &= \int_{\R}q^{-1}   (Jf_1'+Bf_1)^*f_2

%= \int_{\R} q^{-1} f_2^*\ch f_2\ge 0.

%   \end{split}

%\end{equation}
From Lemma \plref{SDM-S1.9} and \eqref{SDM-G1.38c} we infer $\|f_2\|_{A}=0$.
Hence $A f_2=0$ a.e. Since $\cs(J,B,\ch)$
is definite we infer from \eqref{SDM-G1.38a} and \eqref{SDM-G1.38c}
that $f_1=0$. In view of \eqref{SDM-G1.38b} and $A f_2=0$ a.e.
we may apply the definiteness of $\cs(J,B,A)$ to conclude
that $f_2=0$.
\end{proof}

%\newpage
\section[Operators on complete Riemannian manifolds]{First and second order operators on complete Riemannian manifolds}

Let $M$ be a connected complete Riemannian manifold. Furthermore, let
$E$ be a hermitian vector bundle over $M$. We denote by $L^2(E)$
the Hilbert space of square integrable sections of $E$
with respect to the scalar product
\newcommand{\hilbert}[2]{(#1,#2)}
\begin{equation}
      (u,v)=\int_M \scalar{u(p)}{v(p)}_{E_p} d\vol(p).
\label{G1}
\end{equation}
Note that \eqref{G1} is well--defined also if $u$ is
only locally square integrable and $v$ has compact support,
or vice versa. $L^2_{\loc}(E), L^2_{\comp}(E)$ denote
the space of sections of $E$ which are
locally square integrable resp. square integrable with compact 
support. Sometimes it will be convenient to consider distributional
sections of $E$. We denote by $C^{-\infty}(E)$ the (anti)dual
space of $\cinfz{E}$ with respect to the anti-dual pairing
\eqref{G1}. 

Next we consider a second hermitian vector bundle, $F$,
and a first order differential operator
\begin{equation}
   D:\cinfz{E}\longrightarrow \cinfz{F}.
\label{G2}
\end{equation}
Note that we do \emph{not} assume $D$ to be elliptic.
We denote by $D^t$ the formal adjoint of $D$, i.e.
for compactly supported sections $u\in\cinfz{E},v\in\cinfz{F}$
one has
\begin{equation}
    (D u,v)=(v,D^t u).
 \label{G3}
\end{equation}
Thus $D, D^t$ extend to maps on distributional sections of $E,F$
and we will write $Du, D^t v$ also if $u, v$ are distributional
sections of $E, F$, resp. (mostly $u,v$ will at least be locally
square integrable).

Furthermore, let $\hat D$ be the principal symbol of $D$. Then
for $u\in C^{-\infty}(E)$ and $\phi\in\cinf{M}$ one has
\begin{equation}
      D(\phi u)=\hat D(d\phi) u+\phi Du.
      \label{G4}
\end{equation}
\begin{remark}\label{S0}
\eqref{G4} holds whenever all ingredients make sense, 
in particular if $u\in L^2_{\loc}(E),
Du\in L^2_{\loc}(E)$ and $\phi$ is a locally Lipschitz function.
\end{remark}

Note that the defining relation \eqref{G4} for the
principal symbol implies that
\begin{equation}
       \hat D^t(\xi)=-(\hat D(\xi))^*,\qquad \xi\in T_p^*M.
       \label{G5}
\end{equation}

We consider $D$ as an unbounded operator from $L^2(E)$ into
$L^2(F)$. We denote by $\Dmin$ the closure of $D$ and by
$\Dmax=(D^t)^*=((D^t)_{\min})^*$. In general one has
$\Dmin\subsetneqq \Dmax$. Actually, $\Dmin=\Dmax$ is
equivalent to the essential self--adjointness of the
operator
\begin{equation}
     \begin{pmatrix} 0 &D^t\\ D & 0\end{pmatrix}.
\label{G6}
\end{equation}

Next we consider the Schr\"odinger operator
\begin{equation}
     H:=D^tD+V,\label{G9}
\end{equation}
where $V\in L^\infty_{\loc}(\End(E))$ is a locally bounded self--adjoint
(i.e. for each $p\in M$ the endomorphism $V(p):E_p\to E_p$ is self--adjoint)
potential.

$H$ is a symmetric operator in $L^2(E)$ with domain $\cinfz{E}$.
As for $D$ we denote by $H_{\min}$ the closure of $H$ and
$H_{\max}=H^*=H_{\min}^*$.

\begin{dfn}\label{S3} Let $M$ be a complete Riemannian manifold
and let $0<\varrho\le 1$ be a locally Lipschitz function.
We write
\begin{equation}
   \int^\infty \varrho ds=\infty,
\end{equation}
if $\int_0^\infty \varrho(\gamma(t))|\gamma'(t)|dt=\infty$ for any parametrized curve
$\gamma:[0,\infty)\to M$ satisfying $\lim\limits_{t\to\infty} \gamma(t)=\infty$.
The latter limit is taken in the one--point compactification of
$M$, i.e. $\gamma(t)$ eventually leaves any compact subset $K\subset M$.
\end{dfn}

Finally, put $c(x):=\max(1,|\hat D(x)|)$. $c(x)$ is an upper estimate for the
propagation speed of $D$.
Now we can state the main result of this section:

\pagebreak[3]

\begin{theorem}\label{main} Let $q\ge 1$ be a locally Lipschitz function
such that $V\ge -q$. Moreover, assume that 
\begin{thmenum}
\item $c|d(q^{-1/2})|\le C$,
\item $\displaystyle \int^\infty \frac{ds}{c\sqrt{q}}=\infty$,
\item if $u\in\cd(H_{\max})$ then $Du\in L^2_{\loc}(F)$.
\end{thmenum}
Then the operator $H$ is essentially self--adjoint on $\cinfz{E}$.
\end{theorem}

We comment on the assumptions and discuss some special cases:

\begin{remark}1. We emphasize, that the method presented here
is essentially the one of Shubin \cite{Shu:CQC,Shu:ESA},
modulo necessary changes due to the more general class
of operators under consideration.
We found it however worthwhile to show that in principle
all operators of the form $D^tD+V$ can be dealt with in
a unified way, going much beyond the class of Laplace type
operators.

Note also the similarity between Theorem \plref{main} and
Theorem \plref{SDM-S1.10}. Theorem \plref{SDM-S1.10}, in fact,
was inspired by Theorem \plref{main}.

2. The assumption (3) is automatically fulfilled if $D^tD$ is elliptic, or,
more generally, if $D^tD$ is elliptic on a "sufficiently large" subset 
(see Proposition \plref{S2} below). We tried hard to prove the following conjecture:

\begin{conjecture}\label{mainconj} Let $T:\cinfz{E}\to\cinfz{E}$ be a first order differential operator on a
Riemannian manifold and assume that $T^2$ is essentially self--adjoint.
Let $u\in L^2_{\loc}(E), T^2u\in L^2_{\loc}(E)$. Then $Tu\in L^2_{\loc}(E)$.
\end{conjecture}

Let us first comment on why this conjecture is conceivable. If $T^2$
is essentially self--adjoint then $T$ is also essentially self--adjoint
and $\ovl{T^2}=\ovl{T}^2$. Hence, if $u\in L^2(E), T^2 u\in L^2(E)$ then
\begin{equation}\begin{split}
u&\in\cd(\ovl{T^2})=\bigsetdef{v\in L^2(E)}{T^2v\in L^2(E)}\\
                  &=\cd(\ovl{T}^2)=\bigsetdef{v\in L^2(E)}{
Tv, T^2v\in L^2(E)}.\end{split}
\end{equation}
Consequently, $Tu\in L^2(E)$. So, if we remove the "loc" subscripts then
the statement of the conjecture holds. Now, since $T$ is a differential operator,
it is hard to believe that the validity of the conclusion depends on
global properties of $u$. If one believes that the statement
is a purely local one then it should be true even without the essential
self--adjointness assumption on $T^2$, since every symmetric first order
differential operator $T$
can be altered outside a compact set in such a way that all powers
become essentially self--adjoint (cf. the proof of Proposition \plref{S2} below).
Maybe it is possible to prove (or disprove) the conjecture by micro-local
methods. This we did not try too hard.
 
In Proposition \plref{S2a} below it is proved that the conjecture in conjunction
with condition (2) implies condition (3). 

3. Let $V=0$ and $q=1$. Then we obtain the essential self--adjointness
of $D^tD$ if $\int^\infty \frac 1c=\infty$. This is exactly Chernoff's
condition \cite[Thm. 1.3]{Che:ESP}. Note that if $D^tD$ is elliptic then
our method of proof is independent of Chernoff's paper. If $D^tD$ is non--elliptic
we have to use Chernoff's results in the proof of Proposition \plref{S2}
(and also in the proof of Proposition \plref{S2a}).
It is an interesting question whether this Proposition
could be proved by more elementary means.

If $D$ is a generalized Dirac operator then $D$ is elliptic
and $c=1$. Hence we obtain the essential self--adjointness of
$D^2$ (and thus of $D$, too). In this case,
however, our proof is very similar to the one of Wolf
\cite{Wol:ESA}.

4. If $c=1$ then Theorem \plref{main} contains the main results
in \cite{Ole:ESA,Ole:CCQ,Ole:ESAG,Shu:CQC,Shu:ESA,Bra:SAS}
as special cases. Note that loc. cit. mostly deal with cases
where $D^tD$ is a generalized Laplace operator. In this case,
the integrand of $(Hu,v)-(u,Hv)$ can be expressed explicitly
in terms of a divergence. These explicit
divergence formulas are used in an essential way. 
We emphasize that our method works without
such explicit formulas. 
The substitute for them is a more
elaborate use of the calculus of unbounded operators 
in Hilbert space.

In particular, we wanted to include all Dirac type operators. 
For those, of course, the explicit
divergence formulas could be worked out, although it would
be somewhat tedious.

The magnetic Schr\"odinger operator considered in 
\cite{Shu:ESA} is a priori not covered by Theorem \plref{main}
if the magnetic potential is not
smooth. However, if $D^tD$ is elliptic, our proof can easily
be adapted to the case that the $0$th order part of $D$ is
only Lipschitz. For the sake of a simpler presentation,
however, we will confine ourselves to
the case of an operator $D$ with smooth coefficients.
\end{remark}

%\newpage
\subsection{Some Preparations}

\eqref{G3} holds in greater generality:
\begin{lemma}\label{S1} Let $u\in\cd(\Dmax)\cap L^2_{\comp}(E)$ 
and $v\in L^2_{\loc}(F)$ such that $D^tv\in L^2_{\loc}(F)$.
Then $u\in\cd(\Dmin)$ and
\begin{equation}
        (Du,v)=(u,D^tv).\label{G7}
\end{equation}
\end{lemma}
\begin{proof} $u\in\cd(\Dmin)$ follows easily by means of
a Friedrich's mollifier constructed in a neighborhood of
the compact support of $u$.

Next choose a cut--off function $\phi\in\cinfz{M}$ with
$\phi\equiv 1$ in a neighborhood of $\supp u$. Then,
$\phi v\in\cd(D^t_{\max})$ and hence
\begin{equation}\begin{split}
      (Du,v)&= (\phi Du,v)\\
            &= (\Dmin u,\phi v)\\
            &= (u, D^t_{\max} \phi v)\\
            &= (u, -\hat D(d\phi)^* v+ \phi D^t v)\\
            &= (u,D^t v),
		\end{split}
\label{G8}
\end{equation}
since $\supp d\phi\cap \supp u=\emptyset$.
\end{proof}

\begin{lemma}[{cf. Lemma \ref{SDM-S1.8}}]\label{S4} Let $\varrho\ge 1$ be a locally Lipschitz
function on $M$ with $\int^\infty \frac{ds}{\varrho}=\infty$.
Then there is a sequence of Lipschitz functions $(\phi_n)$
with compact support satisfying
\begin{equation}%\begin{split}
         0\le \phi_n\le 1,\quad |d\phi_n|\le \frac{1}{\varrho n},\quad
         \lim_{n\to\infty} \phi_n(x)=1,\quad x\in M.
		%\end{split}
\label{G3.1}
\end{equation}
\end{lemma}
\begin{proof}
Denote by $d_{\varrho}$ the distance function with respect to the
metric $g_\varrho=\varrho^{-2} g$. Then fix $x_0\in M$ and
put $P(x)=d_{\varrho}(x,x_0)$.
As in \cite{Shu:CQC} one concludes $\lim\limits_{x\to\infty}P(x)=\infty$
and $|dP|\le \varrho^{-1}$.
Now choose a cut--off function $\chi\in\cinfz{\R}$ with $0\le \chi\le 1,$
$\chi=1$ near $0$, and $|\chi'|\le 1$. Then put
\begin{equation}
    \phi_n(x)=\chi(\frac{P(x)}{n}).
\end{equation}
$\phi_n$ obviously has the desired properties.
\end{proof}

\begin{prop}\label{S5}Assume that $\int^\infty \frac{ds}{c}=\infty$
and $c|d(q^{-1/2})|\le C$. Let $u\in\cd(H_{\max})$ and $Du\in L^2_{\loc}(F)$. 
Then we have $q^{-1/2}Du\in L^2(F)$ and
\begin{equation}
            \|q^{-1/2}Du\|\le 2\Bigl((1+2C^2) \|u\|^2+ \|u\| \|Hu\|\Bigr).
\end{equation}
\end{prop}
\begin{proof}
Let $0\le \psi\le q^{-1/2}$ be a locally Lipschitz function
with compact support and put $\tilde C=\sup_{p\in M}c(p)|d\psi(p)|$.

%In view of Proposition \plref{S2} we have $Du\in L^2_{\loc}(E)$ and 
Using Lemma \plref{S1} we find
\begin{equation}\begin{split}
     (\psi Du, \psi Du)& = (D^t \psi^2 Du,u)\\
                      &=  2(\psi \hat D^t(d\psi)Du,u)+(\psi^2 D^tDu,u)\\
                  &= 2(\psi \hat D^t(d\psi)Du,u)+(\psi Hu,u)-(V\psi u,\psi u)\\
                  &\le 2\tilde C \|u\| \|\psi Du\| + \|u\| \|Hu\| + \|\psi
                      q^{1/2} u\|^2\\
                  &\le 2\tilde C \|u\| \|\psi Du\| + \|u\| \|Hu\| + \|u\|^2.
		\end{split}
\end{equation}
Using $2 |ab|\le {a^2+b^2}$ the latter can be estimated
\begin{equation}
      \|\psi Du\|^2\le (1+2\tilde C^2) \|u\|^2+\frac 12 \|\psi Du\|^2+ \|u\|
      \|Hu\|,
\end{equation}
and thus
\begin{equation}
      \|\psi Du\|^2\le 2\Bigl((1+2\tilde C^2) \|u\|^2+ \|u\|\|Hu\|\Bigr).
\end{equation}

We apply Lemma \plref{S4} with $\varrho=c$ and obtain a sequence
$(\phi_n)$ of Lipschitz functions $\phi_n$ which satisfy
\eqref{G3.1} with $\varrho=c$. Putting $\psi_n=\phi_n q^{-1/2}$
we have $0\le\psi_n\le q^{-1/2}$ and
\begin{equation}\begin{split}
    c|d\psi_n|&\le cq^{-1/2} |d\phi_n|+ \phi_n c|d(q^{-1/2})|\\
              &\le \frac 1n + C.
		\end{split}
\end{equation}
Since $\psi_n(p)\to q^{-1/2}(p)$ as $n\to\infty$ we reach the conclusion
by invoking the dominated convergence theorem.
\end{proof}

\subsection{Proof of the Main Theorem \ref{main}}
Let $u,v\in\cd(H_{\max})$ and let $0\le \phi$ be a 
Lipschitz function with compact support. Since $q\ge 1$ the condition (2)
implies for any curve $\gamma:[0,\infty)$ as in Definition \plref{S3}
\begin{equation}
\int_0^\infty \frac{1}{c(\gamma(t))}|\gamma'(t)|dt\ge
\int_0^\infty \frac{1}{c(\gamma(t))\sqrt{q(t)}}|\gamma'(t)|dt=\infty,
\label{G3.11}
\end{equation}
hence we can apply Proposition \plref{S5} and find that
$q^{-1/2}Du, q^{-1/2}Dv\in L^2(F)$.
%, in particular $Du,Dv\in L^2_{\loc}(F)$. 
Moreover, since $\phi$ has compact support,
we have $\hat D(d\phi)u\in L^2_{\comp}(F)$.
Also, since $V$ is locally bounded, $D^tD u,D^tD v\in L^2_{\loc}(E).$ 
Finally, the latter implies in view of 
\begin{equation}
      D^t\phi Du=-\hat D(d\phi)^* Du+ \phi D^tDu\in L^2(E).
\end{equation}
Using Lemma \plref{S1} and Remark \plref{S0} we calculate
\begin{equation}\begin{split}
       (\phi u,D^tD v)&=(D\phi u,Dv)\\
                         &=(\hat D(d\phi)u,Dv)+(\phi Du,Dv),
                 \end{split}
\end{equation}
and, similarly,
\begin{equation}
       (D^tDu,\phi v)=(Du,\hat D(d\phi)v)+(\phi Du,Dv).
\end{equation}
Taking differences we obtain
\begin{equation}\begin{split}
      |(\phi u, Hv)-&(Hu,\phi v)|\le |(\hat D(d\phi)u,Dv)|+|(Du,\hat
       D(d\phi)v)|\\
                &\le \sup_{p\in M}\Bigl(|q^{1/2}(p) |\hat D(d\phi)|\Bigr)
       \bigl(\|u\|\|q^{-1/2}Dv\|+\|q^{-1/2}Du\|\|v\|\bigr).
		\end{split}
\end{equation}
Finally we invoke Lemma \plref{S4} with $\varrho=cq^{1/2}$
and choose a sequence of Lipschitz functions $\phi_n$ with compact support
satisfying $0\le \phi_n\le 1$, $|d\phi_n|\le \frac{1}{n c \sqrt{q}}, 
\lim\limits_{n\to\infty} \phi_n(p)=1, p\in M$. Then by dominated convergence
we have on the one hand
\begin{equation}
     (\phi_n u, Hv)-(Hu,\phi_n v)\longrightarrow  (u, Hv)-(Hu,v),\qquad
     n\to\infty,
\end{equation}
and on the other hand
\begin{equation}
|(\phi_n u, Hv)-(Hu,\phi_n v)|\le \frac 1n\bigl( \|u\|\|q^{-1/2}Dv\|+\|q^{-1/2}Du\|\|v\|\bigr).
\end{equation}
This proves the claim.

\subsection{On condition \textup{(3)} and Conjecture \ref{mainconj}}

\begin{prop}\label{S2} Assume that there are compact subsets $K_n\subset M$
such that 
\begin{thmenum}
\item $K_n\subset K_{n+1}$,
\item $\bigcup\limits_{n=1}^\infty K_n=M$,
\item there is an open neighborhood $U_n\supset K_n$ such that 
$D^tD$ is elliptic in $U_n\setminus K_n$.
\end{thmenum}
Let $u\in\cd(H_{\max})$. Then $Du\in L^2_{\loc}(F)$.
\end{prop}
\begin{proof}
1. We note first that if $D^tD$ is elliptic (everywhere)
then this is an easy consequence of elliptic regularity. Namely,
if $Hu=v\in L^2_{\loc}(E)$ then $D^tDu=v-Vu\in L^2_{\loc}(E)$
and hence by elliptic regularity this implies $u\in H^2_{\loc}(E)$. I.e.
$u$ is locally of Sobolev class $H^2$ and hence in particular
$Du\in L^2_{\loc}(E).$

2. If $D^tD$ is not elliptic everywhere then we have to invoke the hyperbolic equation
method as presented e.g. by P. R. Chernoff \cite{Che:ESP}. 
As in 1. we have $D^tDu\in L^2_{\loc}(E)$
and hence, by elliptic regularity, $u\restr U_n\setminus K_n$
is locally of Sobolev class $H^2$, in particular $(Du)\restr U_n\setminus K_n$
is locally square integrable.

We now show that $(Du)\restr K_n$ is square integrable. Choose
a large compact set $K\supset U_n$ and
let $\tilde D$ be a first order differential operator
which coincides with $D$ over $K$ and which vanishes outside a large compact
set $L$. Now consider the operator
\begin{equation}
  T:=\begin{pmatrix} 0 &\tilde D^t\\ \tilde D & 0\end{pmatrix}.
\label{G10}
\end{equation}
$T$ is a formally symmetric differential operator which vanishes outside
a compact set. Hence, $T$ has bounded propagation speed, in particular it
satisfies Chernoff's condition $\int^\infty\frac{ds}{c}=\infty$.
Thus by the hyperbolic equation method \cite{Che:ESP} all
powers of $T$ are essentially self--adjoint. 

Next choose a cut--off function $\phi\in\cinfz{M}$ with $\phi\equiv 1$
in a neighborhood of $K_n$ and $\supp \phi\subset U_n$. Then
the commutator $[D^tD,\phi]=[\tilde D^t\tilde D,\phi]$ is a first 
order differential operator which is supported in $U_n\setminus K_n$.
In particular $[\tilde D^t\tilde D,\phi]u\in H^1_{\comp}(E)$ and
hence $\tilde D^t\tilde D(\phi u)=[\tilde D^t\tilde D,\phi]u+\phi
\tilde D^t\tilde Du=[\tilde D^t\tilde D,\phi]u+\phi D^tDu\in L^2(E)$.
Then
\begin{equation}
      T^2\begin{pmatrix} \phi u \\ 0\end{pmatrix}
          =\begin{pmatrix} \tilde D^t\tilde D \phi u \\0\end{pmatrix}
\end{equation}
is square integrable. Since $T^2$ is essentially self--adjoint, 
this implies that 
\begin{equation}\begin{pmatrix} \phi u \\ 0\end{pmatrix}
\in\cd(\ovl{T^2})=\cd(\ovl{T}^2)=\bigsetdef{v\in L^2(E\oplus E)}{
Tv, T^2v\in L^2(E\oplus E)},
\end{equation}
hence
\begin{equation}
T\begin{pmatrix} \phi u \\0 \end{pmatrix}=\begin{pmatrix} 0\\D(\phi u) \end{pmatrix}
\end{equation}
is square integrable. This implies that $(Du)\restr K_n$ is square integrable.
\end{proof}

\begin{remark} If $D^tD$ is not elliptic in the shells
$U_n\setminus K_n$ then in the proof of 2. we face
the difficulty that there is no obvious way to construct
enough cut--off functions $\phi$ such that
$D^tD(\phi u)\in L^2(E)$. It would be enough
to show the following: given $u\in L^2_{\loc}(E), D^tDu\in L^2_{\loc}(E)$
then there is a $v\in L^2_{\comp}(E),D^tDv\in L^2_{\comp}(E)$ such
that $v\restr K_n =u$. $v$ does not necessarily have to be of
the form $\phi u$.
\end{remark}

\begin{prop}\label{S2a} Assume that Conjecture \plref{mainconj} holds.
Then condition \textup{(2)} in Theorem \plref{main} implies condition
\textup{(3)}.
\end{prop}
\begin{proof} Since $q\ge 1$ the condition (2)
implies $\int^\infty \frac{ds}{c}=\infty$ (cf. \eqref{G3.11}), 
hence the symmetric operator \eqref{G10} satisfies Chernoff's condition
\cite[Thm. 1.3]{Che:ESP}. Thus all powers of $T$ are essentially self--adjoint.
Now, if $u\in \cd(H_{\max})$ then $D^tDu\in L^2_{loc}(E)$ and hence
\begin{equation}
    \tilde u:=\begin{pmatrix} u\\ 0 \end{pmatrix}
\end{equation}
satisfies $\tilde u\in L^2_{\loc}(E\oplus E), T^2 \tilde u\in
L^2_{\loc}(E\oplus E)$.
Consequently, Conjecture \plref{mainconj} implies 
\begin{equation}
    T\tilde u=\begin{pmatrix} 0\\ D u\end{pmatrix}\in L^2_{\loc}(F\oplus F),
\end{equation}
and thus $Du\in L^2_{\loc}(F)$.\end{proof}

%\bibliography{mlabbr,books99,papers99}
%\bibliographystyle{lesch}

\end{document}